\documentclass[draft, 12pt]{amsart}
\usepackage[english]{babel}

\usepackage{amssymb}
\usepackage{amsmath}
\usepackage[all,cmtip]{xy}

\usepackage{color}
\newtheorem{thm}{Theorem}
\newtheorem{defi}{Definition}
\newtheorem{prop}[thm]{Proposition}

\newtheorem{lem}[thm]{Lemma}
\newenvironment{pf}[1]{\noindent\textit{Proof.} #1}{\hfill$\Box$ \medskip}
\newenvironment{rema}[1]{\noindent {\em Remark.} #1}{}

\usepackage{ifdraft, blindtext, verbatim}
\ifdraft{%
    \usepackage{todonotes}}{ 
    \usepackage[disable]{todonotes} 
}

\begin{document}
\title{Hamiltonian loops on the symplectic blow up along a submanifold
}


\author{Andr\'es Pedroza}
\address{Facultad de Ciencias\\
           Universidad de Colima\\
           Bernal D\'{\i}az del Castillo No. 340\\
           Colima, Col., Mexico 28045}
\email{andres\_pedroza@ucol.mx}

\begin{abstract}
We prove that the fundamental group of the group of Hamiltonian diffeomorphisms
of the symplectic manifold that is obtain by blowing up a submanifold  contains an element of infinite order.
We prove this using Weinstein's morphism and 
by constructing explicitly such loop of Hamiltonian diffeomorphisms.
\end{abstract}


\thanks{The author was supported by CONACYT-CB-2017-2018-A1-S-8830 grant}


\maketitle

\section{Introduction}
\label{s:intro}

In \cite{pea-rankham}, it was shown that one-point blow up $(\widetilde M,\widetilde\omega_\rho)$ of a 
closed 4-dimensional symplectic manifold $(M,\omega)$ has a loop $\{\widetilde \psi_t\}$ of Hamiltonian
diffeomorphisms such that the class $[\widetilde\psi_t]$  has infinite order in $\pi_1(\textup{Ham}  
(\widetilde M,\widetilde\omega_\rho)).$  Actually the argument presented in \cite{pea-rankham}
works for any dimension as long as it is greater than two. The restriction to the 4-dimensional case had a 
specific purpose;  it was  related to the  study of the abelian group 
$\pi_1(\textup{Ham} (\widetilde M,\widetilde\omega_\rho)).$ 
Namely, for any positive integer $k$
that there exists closed symplectic 4-manifold such that the rank of
 $\pi_1(\textup{Ham} )$ is at least $k$. 
This phenomena about  the  rank of $\pi_1(\textup{Ham} (\widetilde M,\widetilde\omega_\rho))$ being
positive was first discovered by D. McDuff
in \cite{mcduff-blowup}; the difference between the results of  \cite{pea-rankham} 
and those of \cite{mcduff-blowup}
is that the element of infinite order in $\pi_1(\textup{Ham})$  is written explicitly
in \cite{pea-rankham}. 
The goal of this article consists on 
extending this result. Instead of blowing up a zero-dimensional
submanifold, we blow up a symplectic submanifold $N\subset (M,\omega)$ of codimension $2k>2$.
In this way we also obtain a loop of Hamiltonian diffeomorphisms on $(M,\omega)$, that induces an element
of infinite order in $\pi_1(\textup{Ham}(\widetilde M,\widetilde\omega_{N,\rho}))$
where $(\widetilde M,\widetilde\omega_{N,\rho})$ stands for the symplectic blow up
along the submanifold $N$.

\begin{thm}
\label{t:mainrank}
Let $(M,\omega)$ be a rational closed symplectic manifold and $N\subset (M,\omega)$
a closed symplectic submanifold of codimension $2k>2$.
Let $(\widetilde M,\widetilde\omega_{N,\rho})$ be the symplectic blow up
of $(M,\omega)$ along $N$. Then for some small values of $\rho$ the rank of 
$\pi_1(\textup{Ham}(\widetilde M,\widetilde\omega_{N,\rho}))$
is positive.
\end{thm}

In $(\widetilde M,\widetilde\omega_{N,\rho})$, the symplectic blow
up of   $(M,\omega)$ along $N$,
the parameter $\rho$ in the description of $(\widetilde M,\widetilde\omega_{N,\rho})$
is related to the $\omega^n$-volume of a tubular neighborhood of 
$N\subset(M,\omega)$. Therefore, in the case when $N$ is a single point,
it corresponds to the weight of the one-point blow up. 

Although, Theorem  \ref{t:mainrank} is  a statement about the topology of the
group $\textup{Ham}(\widetilde M,\widetilde\omega_{N,\rho})$, that is as much 
as we know about its topology. 
Nonetheless, the study of the homotopy type of $\textup{Ham}(M,\omega)$ 
has succeeded in some cases, all of which are 4-dimensional. 
For instance the homotopy type of 
 $(\mathbb{C}P^2,\omega_{FS})$
and $(\mathbb{C}P^1\times\mathbb{C}P^1,\omega_{FS}\oplus \omega_{FS})$
was computed by M. Gromov \cite{gromov-psudo}; 
 in \cite{abreu-mcduff-topology},
M. Abreu and D. McDuff 
settled the problem for the one-point blow up of  $(\mathbb{C}P^2,\omega_\textup{FS})$.
For other examples of  4-dimensional symplectic manifolds see  the works of
S. Anjos and S. Eden \cite{anjos-eden},
S. Anjos and M. Pinsonnault \cite{anjos-pin-homo}, J. Evans \cite{Evans-symplectic-mapping} and 
F. Lalonde and M. Pinsonnault  \cite{Lalonde-Pin}. 
As for computations of the fundamental group of $\textup{Ham}(M,\omega)$, there are the works of 
J. Li and  T.-J. Li  \cite{lili-symp}  and J. Li, T.-J. Li and W. Wu \cite{liliwu}.

The techniques that are used to prove Theorem \ref{t:mainrank} are soft techniques.
Furthermore, we give an explicit description of the loop $\{\widetilde\psi_{N,t}\}$ of Hamiltonian 
diffeomorphisms on $(\widetilde M,\widetilde\omega_{N,\rho})$ that induces the element of
infinite order in its fundamental group. In fact, the loop $\{\widetilde\psi_{N,t}\}$ is induced from
a Hamiltonian loop $\{ \psi_t\}$  on $(M,\omega)$ that is supported in a tubular 
neighborhood of $N$. Finally, in order determine that  an element on the fundamental group
of $\textup{Ham}(M,\omega)$ has infinite order, we relay on  
Weinstein's morphism, 
$$
\mathcal{A}:\pi_1(\textup{Ham}(M,\omega) ) \to \mathbb{R}/\mathcal{P}(M,\omega).
$$ 
In this direction we have the following relation
between Hamiltonian diffeomorphism on $(M,\omega)$ and on
$(\widetilde M,\widetilde\omega_{N,\rho})$ for a particular
class of Hamiltonian loops. In the formula of the above result,  
${\Phi(\mathcal{U}_{\rho}(\nu_0))}$ stands for
a particular tubular neighborhood of the symplectic submanifold $N\subset (M,\omega)$.

\begin{thm}
\label{t:weinsrelation}
Let $(M,\omega)$ be a closed symplectic manifold and $N\subset (M,\omega)$
a symplectic submanifold of codimension $2k>2$.
If $\{\psi_t\}$ is a loop of Hamiltonian diffeomorphisms on $(M,\omega)$
that is $N$-liftable to a loop $\{\widetilde\psi_{N,t}\}$
on $(\widetilde M_N,\widetilde \omega_{N,\rho})$, then
\begin{eqnarray*}
\mathcal{A}([\widetilde\psi_{N,t}]) = 
\left[\mathcal{A}([\psi_t]) +
\frac{1}{\textup{Vol}(\widetilde M_N,(\widetilde\omega_{N,\rho})^n)}
\int_0^1 
\int_{\Phi(\mathcal{U}_{\rho}(\nu_0))}  H_t\, \omega^n dt
 \right] 
\end{eqnarray*}
in $\mathbb{R}/\mathcal{P}(\widetilde M_N,\widetilde \omega_{N,\rho})$, where
$H_t$ if the normalized Hamiltonian function of the loop $\{\psi_t\}$.
\end{thm}

Loosely speaking, the reason why the results of \cite{pea-rankham} extend
to the case of a submanifold of positive dimension is because many of the 
arguments that are involved in the definition of symplectic blow up are in a
sense $U(k)$-equivariant. 
Furthermore, the diffeomorphisms
$\{\psi_t\}$ and $\{\widetilde\psi_{N,t}\}$ that appear in the above result are induced
by matrices that lie in the center of $U(k)$. This condition turns to be crucial,
since in many instances the commutativity of the matrices is fundamental.

The paper is organized as follows.  In Section \ref{sec:blow} we review the
 definition of the symplectic one-point blow up as well as some facts 
 of \cite{pedroza-hamiltonian} and \cite{pea-rankham}. 
The definition of the symplectic blow up 
 a long a submanifold is review in Section
\ref{s:blowN}. In Section \ref{s:hamlif} we give conditions that guarantee
that a Hamiltonian diffeomorphisms can be lifted to the symplectic blow up.
Finally, Sections \ref{s:wein} and \ref {s:main} deal with normalization
of Hamiltonian functions and the proofs of the main results respectively.
 
\section{The case of a one-point blow up}
\label{sec:blow}

\subsection{Symplectic structure on the one-point blow up.}
In this section we review how the  symplectic structure on the blow up of the
origin in $(\mathbb{C}^k,\omega_0)$ is defined. 
Following  \cite[Sec. 2]{pedroza-hamiltonian},
we also review the conditions that a Hamiltonian diffeomorphism
on $(\mathbb{C}^k,\omega_0)$ needs to satisfy in order to induced 
 a Hamiltonian diffeomorphism on the one-point blow up.
Later, we will impose similar conditions on a Hamiltonian 
diffeomorphism in order to induced a Hamiltonian on  the symplectic manifold
that is obtain by blowing up a compact symplectic submanifold. 

As a manifold, the one-point blow up of $\mathbb{C}^k$  at the origin
 is the total space of the tautological line bundle,
\begin{equation*}
	\widetilde{\mathbb{C}}^k:=\{(z,\ell)\in \mathbb{C}^k\times \mathbb{C}P^{k-1}  \mid z\in \ell\}.
\end{equation*}
The projection maps are denoted by $\textup{pr} \colon \widetilde{\mathbb{C}}^k\to \mathbb{C}P^{k-1}$ and 
$\pi \colon \widetilde{\mathbb{C}}^k\to \mathbb{C}^k$. 
Here, $\pi$ and $E:=\pi^{-1}(0)$ are called
 the blow up map  is  the exceptional divisor respectively.
Let $\omega_\textup{FS}$ be the 
Fubini-Study symplectic form on $\mathbb{C}P^{k-1}$ normalized so that
 $\langle \omega_\textup{FS}, \mathbb{C}P^1 \rangle = \pi$. 
For each $\rho>0$ define on $\widetilde{\mathbb{C}}^k$     the symplectic form
\begin{equation*}
	\omega_{\rho}:=\pi^*\omega_0+ \rho^2\textup{pr}^*\omega_\textup{FS}.
\end{equation*}
Notice that the $\omega_{\rho}$-area of  a line in the exceptional divisor is $\pi \rho^2$.
Denote by  $B_r\subset \mathbb{C}^k$  the open ball of radius $r$ centered at the origin 
and  its preimage under $\pi$  by
$$
	L_r := \pi^{-1}(B_r).
$$

According to \cite[Prop. 7.1.13]{MS-book} for each $\epsilon\in (0,1)$ there exists 
a symplectic form $\widetilde{\omega}_{\rho,\epsilon}$ on $\widetilde{\mathbb{C}}^k$ such that
\begin{itemize}
\item $\widetilde{\omega}_{\rho,\epsilon}=\omega_\rho$ on $L_\epsilon$ and
\item $\widetilde{\omega}_{\rho,\epsilon}=\pi^{*}(\omega_0)$ on $\widetilde{\mathbb{C}}^k \setminus
\overline{L}_{\rho+\epsilon}.$
\end{itemize}
Later in this note it will be important to know specific details on how  the symplectic form 
$\widetilde\omega_{\rho,\epsilon}$ is defined. To that end, one considers
a diffeomorphism 
$F_{\rho,\epsilon}: \mathbb{C}^k\setminus \{0\}\to \mathbb{C}^k\setminus \overline{B}_\rho$
that consists of stretching the punctured space,
defined as
\begin{eqnarray}
\label{e:strech}
F_{\rho,\epsilon}(z):=f_{\rho,\epsilon}(|z|)\frac{z}{|z|}
\end{eqnarray}
where  $f_{\rho,\epsilon}:\mathbb{R}\to \mathbb{R}$ is a smooth function
 such that 
$f_{\rho,\epsilon}(x)=\sqrt{\lambda^2+x^2}$ on $(0,\epsilon)$ and
is equal to the identity on 
$[\rho+\epsilon,\infty)$.
Then the above symplectic form 
$\widetilde\omega_{\rho,\epsilon}$ is defined as
$\widetilde\omega_{\rho,\epsilon}:=(\pi^*\circ F_{\rho,\epsilon}^*)(\omega_0)$
on  $\widetilde{\mathbb{C}}^k \setminus E$.

In this case $(\widetilde{\mathbb{C}}^k, \widetilde{\omega}_{\rho,\epsilon})$ is defined 
as  the symplectic blow up of the origin in $( \mathbb{C}^k,\omega_0)$ of weight $\rho$. 
It is important to note that the blow up map $\pi \colon (\widetilde{\mathbb{C}}^k, 
\widetilde{\omega}_{\rho,\epsilon})\to ( \mathbb{C}^k,\omega_0)$
is a symplectic diffeomorphism on the complement of 
$\overline{L}_{\rho+\epsilon} \subset  \widetilde{\mathbb{C}}^k.$

The complete details of symplectic one-point blow up, including the exact  
definitions of $\widetilde \omega_{\rho,\epsilon}$ and $f_{\rho,\epsilon}$
appear in \cite[Ch. 7]{MS-book}.

\subsection{Induced Hamiltonian on the one-point blow up.}
\label{ss:inducedHam}
Next we review some results of \cite{pedroza-hamiltonian} that deal with the problem of lifting a
Hamiltonian diffeomorphism on $( \mathbb{C}^k,\omega_0)$ to a Hamiltonian on
$(\widetilde{\mathbb{C}}^k, \widetilde{\omega}_{\rho,\epsilon})$.

Let $\phi$ and $\widetilde \phi$   be  Hamiltonian diffeomorphisms on $(\mathbb{C}^k,\omega_0)$ and
$(\widetilde{\mathbb{C}}^k, \widetilde{\omega}_{\rho,\epsilon})$ respectively.
We say that the diffeomorphism $\widetilde \phi$ 
{\em lifts} the diffeomorphism $\phi$ if the following diagram commutes
$$
\xymatrix{(\widetilde{\mathbb{C}}^k, \widetilde{\omega}_{\rho,\epsilon})\ar[d]_\pi\ar[r]^{\widetilde{\phi}} &
 (\widetilde{\mathbb{C}}^k, \widetilde{\omega}_{\rho,\epsilon})\ar[d]^\pi \\
({\mathbb{C}}^k, {\omega}_0)\ar[r]^\phi & ({\mathbb{C}}^k, {\omega}_0).
}$$
This means that $	\pi\circ\widetilde{\phi}=\phi\circ\pi$.

Consider $\{\psi_t\}_{0\leq t\leq 1}$ a path of Hamiltonian diffeomorphisms on $(\mathbb{C}^k,\omega_0)$
with Hamiltonian function $H_t$ such that $\psi_0=1_{\mathbb{C}^k}$. 
The goal is to define a path of Hamiltonian diffeomorphisms
$\{\widetilde{\psi}_t\}$
on 	$(\widetilde{\mathbb{C}}^k, \widetilde \omega_{\rho,\epsilon})$ that lifts
the path $\{\psi_t\}$; that is
 $	\pi\circ\widetilde{\psi}_t=\psi_t\circ\pi$
  for each 
$t\in[0,1]$. 
At this point it is important to note that 
the Hamiltonian path induced by the function 
$H_t\circ \pi:(\widetilde{\mathbb{C}}^k, \widetilde{\omega}_{\rho,\epsilon})\to\mathbb{R}$ fails to
lift the path $\{\psi_t\}$. 
Away from  the exceptional divisor is true that the path induced by 
$H_t\circ \pi$ is a lift of the original path; the problem occurs  near the 
exceptional divisor. Note for instance that $H_t\circ \pi$ is independent of the parameter $\rho$
and such parameter is relevant in the symplectic structure on the one-point blow up.
As we will see, the Hamiltonian function that induces the desired Hamiltonian path
will depend  on the parameter $\rho$.

The required condition for a path of Hamiltonian diffeomorphisms to lift to the one-point blow up is 
that for each $t\in[0,1]$
there exists a matrix $U_t\in U(k)$  such that $\psi_t=U_t$ on 
${B}_{\rho+\epsilon}\subset\mathbb{C}^k$.

\begin{lem}
\label{l:hamfunc}
Let $\psi_t :(B_{r},\omega_0)\to (B_{r},\omega_0)$ 
be a Hamiltonian path given by unitary matrices
starting at the identity matrix with Hamiltonian function $H_t$. Then
$$
H_t(z)=H_t(\lambda z)
$$ 
for $z\in B_r$ and $\lambda \in S^1$. 
\end{lem}
\begin{pf}	
Denote by $X_t$ the time-dependent vector field of
the path $\{\psi_t\}$. For $\lambda\in S^1$, let $\phi_\lambda:B_r\to B_r$ be matrix multiplication
by $\lambda I$. Since $\phi_\lambda$ is in the center of $U(k)$ it follows that
\begin{eqnarray*}
X_t\circ \phi_\lambda =\left.\frac{d}{ds}\right|_{s=t} \psi_s\circ \psi_t ^{-1}\circ  \phi_\lambda
=\left.\frac{d}{ds}\right|_{s=t} \phi_\lambda \circ \psi_s\circ \psi_t ^{-1} = (\phi_\lambda)_* X_t.
\end{eqnarray*}

Therefore
$$
d (H_t\circ \phi_\lambda)=\omega_0(X_t, (\phi_\lambda)_*( \cdot)) =
(\phi_\lambda) ^*\omega_0 (X_t,\cdot)=d H_t.
$$
Since  both functions   $H_t$ and $H_t\circ \phi_\lambda$ agree at the origin, it follows that
 $H_t(z)=H_t(\lambda z)$.
\end{pf}

Let $\{\psi_t\}$ be a path of Hamiltonian diffeomorphisms on $(\mathbb{C}^k,\omega_0)$
such that for each $t$
the restriction of $\psi_t$  to ${B}_{\rho+\epsilon}\subset\mathbb{C}^k$
is given by a unitary matrix. If $H_t: (\mathbb{C}^k,\omega_0)\to \mathbb{R}$ is the Hamiltonian 
function of the path, define the function 
$\widetilde H_t:(\widetilde{\mathbb{C}}^k, \widetilde \omega_{\rho,\epsilon})\to \mathbb{R}$
by
\begin{eqnarray}
\label{e:hamblow}
 \widetilde H_t(z,\ell) :=
\left\{
	\begin{array}{ll}
			H_t\circ F_{\rho,\epsilon}\circ \pi (z,\ell) 
			& \textup{if } (z,\ell)\in \widetilde{\mathbb{C}}^k\setminus E \\
			H_t\left(  \frac{\rho}{|w|} w\right) & \textup{if $z=0$ and $[w]=\ell  \in \mathbb{C}P^{k-1}$}.\\
	\end{array}
\right.
\end{eqnarray}
It follows from Lemma \ref{l:hamfunc} and the definition of $F_{\rho,\epsilon}$ that $ \widetilde H_t$
is well defined. Note that $ \widetilde H_t$ depends on the diffeomorphism $F_{\rho,\epsilon}$ and hence
on the parameter $\rho$.
The path of Hamiltonian diffeomorphisms on $(\widetilde{\mathbb{C}}^k, \widetilde \omega_{\rho,\epsilon})$
induced by $\widetilde H_t$
is the path that lifts the initial path on $({\mathbb{C}}^k, {\omega}_0)$.
The proof of this fact  as well as the proof of the next proposition appears in 
\cite[Sec. 3]{pedroza-hamiltonian}.

\begin{prop}
The time-dependent vector field $\widetilde X_t$ and path of diffeomorphisms
$\{\widetilde \psi_t\}$ induced by the function 
$\widetilde H_t:(\widetilde{\mathbb{C}}^k, \widetilde \omega_{\rho,\epsilon})\to \mathbb{R}$,  defined in
(\ref{e:hamblow}), are such that
\begin{enumerate}
\item $\pi_*(\widetilde X_t)=X_t$ for each $t\in[0,1]$ and
\item $\widetilde \psi_t$ is a lift of $\psi_t$ for each $t\in[0,1]$.
\end{enumerate}
In particular, $\{\widetilde \psi_t\}$ is a Hamiltonian path of 
$(\widetilde{\mathbb{C}}^k, \widetilde \omega_{\rho,\epsilon})$ that lifts the Hamiltonian path 
$\{\psi_t\}$.
\end{prop}

Thus if $\psi$ is a Hamiltonian diffeomorphism on $(\mathbb{C}^k,\omega_0)$
such that it can be joined to the identity by a path of Hamiltonian diffeomorphisms such that
each one of them is equal to a unitary matrix on the ball ${B}_{\rho+\epsilon}$, 
then $\psi$ admits a
Hamiltonian lift to the one-point blow up $(\widetilde{\mathbb{C}}^k, \widetilde \omega_{\rho,\epsilon})$.

This definition can be extended to
the one-point blow up of weight $\rho$ of
 an arbitrary symplectic manifold $(M,\omega).$ 
For the purpose of this note it is not relevant to review such definition.
For further details  see \cite{pedroza-hamiltonian}.

\section{Review of the blow up along a submanifold}
\label{s:blowN}

Following \cite[Ch. 7]{MS-book},
we review the construction of the symplectic blow up along a compact symplectic submanifold. 
The complete details of this construction are relevant for the Hamiltonian 
diffeomorphism that we will defined on the blown up manifold.

\subsection{The blow up along a compact submanifold.}
Let $N\subset (M,\omega)$ be a compact symplectic submanifold of codimension $2k>2$.
Then
$\pi_N:TN^\omega\to N$ corresponds to its normal bundle which is
a symplectic vector bundle. Since 
the symplectic linear group contracts to the  unitary group, 
$\pi_N$ is in fact a complex vector bundle.
For simplicity, write $\nu$ for  $TN^\omega$.
Thus there exists a  Hamiltonian linear action of $U(k)$ on $(\mathbb{R}^{2k},\omega_0)
=(\mathbb{C}^{k},\omega_0)$
and  a principal
$U(k)$-bundle $P\to N$  such that the associated fibration coincides with the normal bundle,
$\nu\simeq P\times_{U(k)}\mathbb{C}^{k}$. 
Let $\mu:(\mathbb{C}^{k},\omega_0)\to \mathfrak{u}(k)^*$ be the 
corresponding moment map
and  fix  a connection 1-form $A\in\Omega^1(P, \mathfrak{u}(k))$ on  $P\to N$.
On
$P\times\mathbb{C}^k$ define the 1-form
$$
\langle \mu,A \rangle_{(p,z)} (u+v):=
\langle \mu(z),A_p(u) \rangle
$$ 
where $ u\in T_pP $ and $ v\in  T_z \mathbb{C}^k$. 
Here $\langle\cdot ,\cdot\rangle$
stands for the pairing between  $\mathfrak{u}(k)^*$ and $\mathfrak{u}(k)$. 
Notice that the 1-form $\langle \mu,A \rangle$ is independent of the vectors that lie in $T\mathbb{C}^k$. 
For $\xi\in \mathfrak{u}(k)$, let $X_\xi$ be the vector
on $\mathbb{C}^k$ induced by the action of $U(k)$.
 Then  on  $P\times\mathbb{C}^{k}$ define  the 2-form
\begin{eqnarray*}
\omega_{A}(u_1+v_1,u_2+v_2 )&:=&\omega_0(v_1+X_{A(u_1)},v_2+X_{A(u_2)})-\\
& & d\langle \mu, A\rangle(u_1+v_1,u_2+v_2 ).
\end{eqnarray*}
So defined, the 2-form $\omega_A$ is closed. Moreover, according to  \cite[Thm. 6.6.3]{MS-book} 
$\omega_A$   is $U(k)$-invariant and horizontal.
Henceforth, it induces a  closed  2-form on $\nu\simeq P\times_{U(k)}\mathbb{C}^{k}$. 
In order to avoid cumbersome notation, we will make no distinction between  
a basic form and the induced form on the orbit space.
Note that on each fibre of $\pi_N:\nu\to N$ the form 
$\omega_A$ restricts  to $\omega_0$.
Finally, define the 2-form
\begin{equation}
\label{e:formonN}
\omega_{\nu,A}:=\omega_A  + \pi_N^*(\omega)
\end{equation}
on $\nu.$ Notice that $\omega_{\nu,A}$  also restricts to $\omega_0$ on each fibre.
Furthermore, if  $\nu_0$ stands for the zero-section of $\nu$ then  $\omega_{\nu,A}$
restricted to $\nu_0$ is equal to $\pi_N^*(\omega)$.
These observations imply the following lemma. First, for $\epsilon >0$  denote by
$\mathcal{U}_{\epsilon}(\nu_0):=
P\times_{U(k)} B_{\epsilon}$ the $\epsilon$-disk subbundle of
$\nu\to N$.

\begin{lem}
\label{l:sympformnormal}
Let $N\subset (M,\omega)$ be a compact symplectic submanifold of codimension $2k$. 
Then there exist
$\epsilon_0>0$ such that the 2-form 
$$
\omega_{\nu,A}=\omega_A  + \pi_N^*(\omega)
$$
is symplectic on the
disk subbundle $\mathcal{U}_{\epsilon_0}(\nu_0)$.  Moreover $\omega_{\nu,A}$
restricts to $\omega_0$ on each fibre and to 
$\pi_N^*(\omega)$ on the zero-section.
\end{lem}

It follows by the symplectic neighborhood theorem there exists a symplectic diffeomorphism
\begin{eqnarray}
\label{e:symplecticdiffeo}
\Phi: (\mathcal{U}_{\epsilon_0}(\nu_0),\omega_{\nu,A}) \to (\mathcal{U}(N), \omega) 
\end{eqnarray}
where $\mathcal{U}(N)\subset (M,\omega)$ is some neighborhood
of $N$.

The outline of the construction of the symplectic blow up of $(M,\omega)$ along the symplectic
submanifold 
$N$ is as follows:
The  neighborhood  $ (\mathcal{U}(N), \omega)$ of $N\subset (M,\omega)$ is identified with 
 the open set $(\mathcal{U}_{\epsilon_0}(\nu_0),\omega_{\nu,A})$ 
 via the symplectic diffeomorphism $\Phi$ and thus the zero-section $\nu_0$ in $\nu$ 
is identified with the submanifold $N$ is blown up. Hence we give a description
of the blow up of the zero-section in $\nu$. 
By the blow up of $\nu_0$  we mean the total space of the associated bundle
that is  obtained by $\widetilde{\mathbb{C}}^{k}$ and $P\to N$ which is  denoted by 
$\widetilde\nu$.
Finally, $\widetilde\nu$ carries a symplectic form $\widetilde\omega_{\widetilde\nu,A,\rho}$
and such model
is glued back to $ (M\setminus N,\omega)$.

The complex-linear action of $U(k)$ on $\mathbb{C}^k$ has a canonical lift
 to $\widetilde{\mathbb{C}}^{k}$
defined as
$$
U.(z,\ell)\mapsto (U(z),U(\ell)) \hspace{1.5cm} \textup{ for }  U\in U(k) \textup{ and }(z,\ell)\in 
\widetilde{\mathbb{C}}^{k}.
$$ 
This $U(k)$-space together with the $U(k)$-principal bundle $P\to N$ 
induce the associated fibration 
$\widetilde\pi_N:\widetilde\nu\to N$ where $\widetilde \nu:=P\times_{U(k)}\widetilde{\mathbb{C}}^{k}$.
Since the blow up map $\pi:\widetilde{\mathbb{C}}^{k}\to\mathbb{C}^{k}$ is $U(k)$-equivariant,
 it induces the map 
$$
\widetilde\pi: P\times_{U(k)}\widetilde{\mathbb{C}}^{k}\to P\times_{U(k)}{\mathbb{C}}^{k}
$$
defined as $\widetilde\pi[p,(z,\ell)]:=[p,\pi(z,\ell)]=[p,z].$ Notice that $\widetilde\pi:\widetilde\nu\to \nu$
corresponds to the blow up of the zero-section $\nu_0$ in $\nu$.
As in the case of the one-point blow up, $\widetilde\pi$ is a diffeomorphism on the complement of
$\widetilde{\pi}^{-1}(\nu_0)\subset \widetilde\nu$.

For any $\epsilon>0$, the ball $B_\epsilon\subset\mathbb{C}^k$ is  $U(k)$-invariant 
and hence  $L_\epsilon\subset\widetilde{\mathbb{C}}^{k}$ is also $U(k)$-invariant.
Write ${\mathcal{U}}_\epsilon(\widetilde{\pi}^{-1}(\nu_0))$ for the subbundle
$P\times_{U(k)}L_\epsilon$. Notice that 
 $\widetilde{\pi}$ maps
${\mathcal{U}}_\epsilon(\widetilde{\pi}^{-1}(\nu_0))$
onto $\mathcal{U}_{\epsilon}(\nu_0)$. Furthermore, $\widetilde \pi$
maps  ${\mathcal{U}}_\epsilon(\widetilde{\pi}^{-1}(\nu_0))\setminus \widetilde{\pi}^{-1}(\nu_0)$
diffeomorphically onto   $\mathcal{U}_{\epsilon}(\nu_0)\setminus \nu_0$. Summarizing, we have the 
following commutative diagram of maps
$$
\xymatrix{{\mathcal{U}}_\epsilon(\widetilde{\pi}^{-1}(\nu_0)) \ar[r]^{\widetilde \pi}\ar@{->}
[dr]_{\widetilde \pi_N} &
\mathcal{U}_\epsilon(\nu_0) \ar[d]^{\pi_N}  \ar@{^(->}[r] \ar@{->}
[dr]^{\Phi}  & \nu \\
&N \ar@{^(->}[r] & \mathcal{U}(N)\ar@{^(->}[r]  & M.   
}$$

We are now ready to define the manifold that results by blowing up  $(M,\omega)$ along the submanifold $N$.
If $\epsilon_0>0$ is as before,
the blow up of $M$ along the compact symplectic submanifold $N$ is defined as
\begin{eqnarray}
\label{e:mblow}
\widetilde M_N:=M\setminus N \ \  
\cup_{\Phi\circ \widetilde\pi } \ \ {\mathcal{U}}_{\epsilon_0}(\widetilde{\pi}^{-1}(\nu_0))
\end{eqnarray}
where
${\mathcal{U}}_{\epsilon_0}(\widetilde{\pi}^{-1}(\nu_0))\setminus 
\widetilde{\pi}^{-1}(\nu_0)$ is identified with $\mathcal{U}(N)\setminus N$
 via the diffeomorphism $\Phi\circ \widetilde\pi$.
 
As in the case of the one-point blow up, there is a projection map 
$\widetilde\pi_{(M,N)}: \widetilde M_N\to M$, also called the blow up map,
defined as
\begin{eqnarray}
\label{e:hamblowmap}
 \widetilde\pi_{(M,N)}(p) :=
\left\{
	\begin{array}{ll}
	p & \textup{if } p \in M\setminus N \\
	\Phi\circ \widetilde\pi(p) & \textup{if } 
	p=[p^\prime, (z,\ell)]\in {\mathcal{U}}_{\epsilon_0}(\widetilde{\pi}^{-1}(\nu_0)).
	\end{array}
\right.
\end{eqnarray}
Note that  $\widetilde\pi_{(M,N)}$ depends on the disk subbundle
 ${\mathcal{U}}_{\epsilon_0}(\nu_0)$. We omit such dependency from the notation of the map.
The preimage of the blown up submanifold $E_N:={\widetilde\pi}_{(M,N)}^{-1}(N)$ is called the exceptional
divisor and ${\widetilde\pi}_{(M,N)}$ restricted to the exceptional divisor, 
$E_N\to N$, is a fibration with fiber $\mathbb{C}P^{k-1}$ that corresponds
to the projectivization of the normal bundle $\nu\to M$.
  
It remains to define a symplectic structure on $\widetilde M_N$.
In the definition of $\widetilde M_N$,
the parameter $\epsilon_0>0$ of the disk subbundle was
irrelevant. 
Its importance will become apparent
in the description of the symplectic structure on it.

\subsection{The symplectic structure on $\widetilde M_N$.} 
\label{sec:sympblow}
Let $\widetilde \mu_\rho: (\widetilde{\mathbb{C}}^{k},\widetilde\omega_{\rho,\epsilon})
\to\mathfrak{u}(k)^*$ be the  moment map of the induced action of $U(k)$
on $(\widetilde{\mathbb{C}}^{k},\widetilde\omega_\rho)$ described above.
Using the same connection 1-form $A$
on the principal bundle $P\to N$, we have the analog
of the form $\omega_A$ defined before, but now on 
$P\times\widetilde{\mathbb{C}}^{k}$.
That is, on $P\times\widetilde{\mathbb{C}}^{k}$ we have the
following 2-form
\begin{eqnarray*}
\widetilde\omega_{A,\rho}(u_1+v_1,u_2+v_2 )&:=&
\widetilde\omega_{\rho,\epsilon}(v_1+X_{A(u_1)},v_2+X_{A(u_2)})-  \\
& & d\langle \widetilde\mu_\rho, A\rangle(u_1+v_1,u_2+v_2 ) 
\end{eqnarray*}
that is closed and  descends to the orbit space
$\widetilde\nu:=P\times_{U(k)} \widetilde{\mathbb{C}}^{k}$. Moreover it
 restricts to  $\widetilde\omega_{\rho,\epsilon}$ on each fibre of $\widetilde\pi_N:\widetilde\nu\to N$.
As before, on  $\widetilde\nu$ we define the 2-form
\begin{equation}
\label{e:symplecformblow}
\widetilde\omega_{\widetilde\nu,A,\rho}:=\widetilde\omega_{A,\rho} +\widetilde\pi_N^*(\omega).
\end{equation}

Recall that $\pi: (\widetilde{\mathbb{C}}^{k},\widetilde\omega_{\rho,\epsilon})\to (\mathbb{C}^{k},\omega_0)$
is  a symplectic diffeomorphism on the complement of $\overline{L}_{\rho+\epsilon}$. 

\begin{lem}
\label{l:symprel}
Let $\rho$ and $\epsilon\in(0,1)$ such that
$\rho+\epsilon<\epsilon_0$. Then under the map
$\widetilde \pi: \widetilde\nu\to \nu$ we have that 
$$
\widetilde\pi^*(\omega_{\nu,A})= \widetilde\omega_{\widetilde\nu,A,\rho}
$$
on 
$P\times_{U(k)}(\widetilde{\mathbb{C}}^k\setminus \overline{L}_{\rho+\epsilon})$.
\end{lem}
\begin{proof}
In order to prove the claim, consider the forms on $P\times {\mathbb{C}}^k$ 
instead of $P\times_{U(k)} {\mathbb{C}}^k$. Similarly for  the forms on 
$P\times_{U(k)} \widetilde{\mathbb{C}}^k$. Hence we must show that 
$$
(1\times\pi)^*(\omega_{\nu,A})= \widetilde\omega_{\widetilde\nu,A,\rho}.
$$

First note that since $\pi_N$ and $\widetilde\pi_N$ are fibrations over $N$ and $\widetilde\pi$
is a bundle map, it follows that $\widetilde\pi_N^*(\omega)=   \widetilde\pi^*
\circ\pi_N^*(\omega)$ on all $\widetilde\nu$. 
$$
\xymatrix{
P\times \widetilde{\mathbb{C}}^k \ar[r]\ar[d]_{1\times \pi} & \widetilde\nu
\ar[dr]^{\widetilde\pi_N} \ar[d]_{\widetilde\pi}& \\
P\times {\mathbb{C}}^k \ar[r]                         &   \nu \ar[r]_{\pi_N} & N 
}$$
Therefore 
$(1\times \pi)^*(   \pi_N^*(\omega)  )=   \widetilde\pi_N^*(\omega)$.
It remains to show that 
$
(1\times\pi)^*(\omega_{A})= \widetilde\omega_{A,\rho}.
$ In this step the parameter $\rho$ will be relevant. Since
 $\pi: (\widetilde{\mathbb{C}}^{k},\widetilde\omega_{\rho,\epsilon})\to (\mathbb{C}^{k},\omega_0)$
is  a symplectic diffeomorphism on the complement of $\overline{L}_{\rho+\epsilon}$ and is $U(k)$-equivariant, 
the moment maps satisfy
$\widetilde\mu_\rho=\mu\circ\pi$    on 
$P\times(\widetilde{\mathbb{C}}^k\setminus \overline{L}_{\rho+\epsilon})$.
This in turn imply  that  
$(1\times\pi)^*\langle \mu,A  \rangle =\langle \widetilde\mu_\rho,A  \rangle$
also on 
$P\times(\widetilde{\mathbb{C}}^k\setminus \overline{L}_{\rho+\epsilon})$.

Once again, using the fact that the blow up map is a symplectic diffeomorphism
on the complement of $\overline{L}_{\rho+\epsilon}$ and the fact that $\pi_*(\widetilde X_{\xi})
=X_{\xi}$ it follows that 
 \begin{eqnarray*}
\omega_0(\pi_*(v_1)+X_{A(u_1)},  \pi_*(v_2)+X_{A(u_2)}  )
&=& \omega_0(\pi_*(v_1)+\pi_*(\widetilde{X}_{A(u_1)}),\\ 
& & \hspace{.6cm}\pi_*(v_2)+\pi_*(\widetilde{X}_{A(u_2)})  )\\
&=& \pi^*(\omega_0)(v_1+\widetilde{X}_{A(u_1)}, v_2+\widetilde{X}_{A(u_2)}  )\\
&=& \omega_{\rho,\epsilon}(v_1+\widetilde{X}_{A(u_1)},  v_2+\widetilde{X}_{A(u_2)}  )
\end{eqnarray*}
where $u_j\in TP$ and $v_j\in T\widetilde{\mathbb{C}}^k$.

All these facts  show that 
$(1\times \pi)^*(\omega_{\nu,A})= \widetilde\omega_{\widetilde\nu,A,\rho}$ on
$P\times(\widetilde{\mathbb{C}}^k\setminus \overline{L}_{\rho+\epsilon})$. Henceforth,
$
\widetilde\pi^*(\omega_{\nu,A})= \widetilde\omega_{\widetilde\nu,A,\rho}
$
on 
$P\times_{U(k)}(\widetilde{\mathbb{C}}^k\setminus \overline{L}_{\rho+\epsilon})$.
\end{proof}

For the rest of the article we fix the parameters $\rho>0$ and $\epsilon\in(0,1)$ such that
$\rho+\epsilon<\epsilon_0$.
Remember that 
$
\Phi: (\mathcal{U}_{\epsilon_0}(\nu_0),\omega_{\nu,A}) \to (\mathcal{U}(N), \omega) 
$ is a symplectic diffeomorphism 
and by Lemma \ref{l:symprel}
that 
$
\widetilde\pi^*(\omega_{\nu,A})= \widetilde\omega_{\widetilde\nu,A,\rho}
$
on $P\times_{U(k)}(\widetilde{\mathbb{C}}^k\setminus \overline{L}_{\rho+\epsilon})$.
Therefore on $P\times_{U(k)}(\widetilde{\mathbb{C}}^k\setminus \overline{L}_{\rho+\epsilon})$ 
we have that $  \widetilde\pi^* \circ \Phi^* (\omega) = \widetilde\omega_{\widetilde\nu,A,\rho}$.
Hence it follows that $\widetilde M_N$ carries a symplectic form
$\widetilde \omega_{N,\rho}$ defined as
\begin{eqnarray}
\label{e:symplcformblow}
\widetilde \omega_{N,\rho}:=
\left\{
	\begin{array}{ll}
	\omega & \textup{on }M\setminus\Phi(\overline{ \mathcal{U}_{\epsilon^\prime}(\nu_0)  })  \\
	\widetilde\omega_{\widetilde\nu,A,\rho}& \textup{on }\mathcal{U}_{\epsilon_0}(\widetilde \pi^{-1}(\nu_0))
	\end{array}
\right.
\end{eqnarray}
where $\epsilon^\prime$  is such that $\rho+\epsilon <\epsilon^\prime<\epsilon_0$.
Note that $\widetilde\pi_{(M,N)}$ is a symplectic diffeomorphism on
$M\setminus\Phi(\overline{ \mathcal{U}_{\epsilon^\prime}(\nu_0)  })$.

\section{Induced Hamiltonians on the blown up manifold}
\label{s:hamlif}

In this section we give conditions  on a Hamiltonian diffeomorphism on
 $(M,\omega)$ 
in order to ensure that it lifts to  a Hamiltonian diffeomorphism
on the blown up manifold $(\widetilde M_N, \widetilde\omega_{N,\rho})$.
As stated in the previous section, a diffeomorphism $\widetilde \phi$
on $\widetilde M_N$ is said to be a lift
of the diffeomorphism $\phi$ on $M$ if  the following
diagram commutes
$$
\xymatrix{
\widetilde{M}_N\ar[d]_{\widetilde\pi_{(M,N)}}\ar[r]^{\widetilde\phi} &
\widetilde{M}_N\ar[d]^{\widetilde{\pi}_{(M,N)}}\\
M\ar[r]^\phi & M.
}$$
By the nature of the blow up construction, these conditions will only take place in a 
neighborhood of the symplectic submanifold $N\subset M$. 
To that end,
 fix a tubular neighborhood $\mathcal{U}(N)$ of $N$  and a symplectic
diffeomorphism
$
\Phi: (\mathcal{U}_{\epsilon_0}(\nu_0),\omega_{\nu,A}) \to (\mathcal{U}(N), \omega) 
$ as before.

\begin{defi}
Let $N\subset (M,\omega)$ be a compact symplectic submanifold
of codimension $2k$.
A symplectic diffeomorphism $\psi: (M,\omega)\to (M,\omega)$
is called {\em $N$-liftable} if 
\begin{itemize}
\item[a)] $\psi(\mathcal{U}(N))\subset \mathcal{U}(N)$  and
\item[b)] the map  $\Phi^{-1}\circ \psi\circ\Phi:   
(\mathcal{U}_{\epsilon_0}(\nu_0),\omega_{\nu,A}) \to
 (\mathcal{U}_{\epsilon_0}(\nu_0),\omega_{\nu,A})$ 
 takes the form
$$
 [p,z]\mapsto [p,U(z)]
 $$
where
$U$ lies in the center of $ U(k)$.
\end{itemize}
\end{defi}

\begin{rema}
\begin{itemize}
\item
The above definition depends on the tubular neighborhood $\mathcal{U}(N)$
and  the diffeomorphism $\Phi$. 
In order to avoid cumbersome notation we omit 
such dependencies from the definition.

\item The reason that the matrix $U$ must lie in the center
of $U(k)$ is to ensure that
 $\widetilde\mu_\rho(U(z))=\widetilde\mu_\rho(z)$,
 where $\widetilde \mu_\rho: (\widetilde{\mathbb{C}}^{k},\widetilde\omega_{\rho,\epsilon})
\to\mathfrak{u}(k)^*$ is the moment map of the induced action 
used in the definition of the symplectic form on  
described in Sec. \ref{sec:sympblow}. 
It will also be useful when proving that the induced diffeomorphisms on 
the blown up manifold by a Hamiltonian is also Hamiltonian.

\item If  $\phi$ is a gauge transformation of $P\to N$, then the map
$
 [p,z]\mapsto [\phi(p),U(z)]
$ is a diffeomorphism of $\mathcal{U}_{\epsilon_0}(\nu_0)$. However,
if $\phi\neq 1_P$, such diffeomorphism
on  $(\mathcal{U}_{\epsilon_0}(\nu_0),\omega_{\nu,A})$
will not be symplectic.
\end{itemize}
\end{rema}

\medskip

If $\psi: (M,\omega)\to (M,\omega)$ is $N$-liftable, denote by $\psi_\nu$  
the diffeomorphism $\Phi^{-1}\circ \psi\circ\Phi$ defined on 
$(\mathcal{U}_{\epsilon_0}(\nu_0),\omega_{\nu,A})$.
Hence   
 $\psi_\nu$ induces the diffeomorphism
$\widetilde\psi_\nu  :
\mathcal{U}_{\epsilon_0}(\widetilde \pi^{-1}(\nu_0)) \to \mathcal{U}_{\epsilon_0}(\widetilde \pi^{-1}(\nu_0)) $
given  by
$$
\widetilde\psi_\nu [p, (z, \ell)]:=[p, (U(z), U(\ell))].
$$
Since $\psi$ is $N$-liftable it  maps the complement of  $\; \mathcal{U}(N)\subset M$ to itself. 
Therefore $\widetilde\psi_\nu$ together  with $\psi$ give rise to 
a well-defined diffeomorphism on the
blown up manifold $\widetilde M_N$ that we denote by  $\widetilde\psi_N$ .
So defined $\widetilde\psi_N$ is a lift of $\psi$, that is 
 $\widetilde\pi_{(M,N)}\circ\widetilde\psi_N= \psi \circ\widetilde\pi_{(M,N)}$.
It remains to show that it is a symplectic diffeomorphism.

\begin{prop}
\label{p:liftsymplectic}
If $\psi: (M,\omega)\to (M,\omega)$ is a symplectic diffeomorphism that is 
$N$-liftable, then $\widetilde\psi_N$ is a symplectic diffeomorphism
of $(\widetilde M_N,\widetilde\omega_{N,\rho}).$
\end{prop}
\begin{proof}
By definition of the symplectic form  $\widetilde\omega_{N,\rho}$ in Eq. 
(\ref{e:symplcformblow})
and of the diffeomorphism $\widetilde\psi_N$, if follows that $\widetilde\psi_N$ is a 
symplectic diffeomorphism on 
$M\setminus\Phi(\overline{ \mathcal{U}_{\epsilon^\prime}(\nu_0)  })$. It only remains to
show that $\widetilde\psi_N$ is a 
symplectic diffeomorphism on 
$ \mathcal{U}_{\epsilon_0}(\widetilde \pi^{-1}(\nu_0))$.

Since $\psi$ is $N$-liftable, on $\mathcal{U}_{\epsilon_0}(\widetilde \pi^{-1}(\nu_0))$ the diffeomorphism
$\widetilde\psi_N$ is given by $ \widetilde \psi_\nu$. That is 
$$
\widetilde \psi_\nu[p,(z,\ell)]=[p,(U(z),U(\ell))]
$$
where $U$ lies in the center of $U(k)$.
In order to show that $\widetilde\psi_\nu^*(\widetilde\omega_{\widetilde\nu,A,\rho})=
\widetilde\omega_{\widetilde\nu,A,\rho}$, we regard  the symplectic form
on $P\times L_{\epsilon_0}$ instead of 
$P\times_{U(k)}L_{\epsilon_0}$. In particular, $\widetilde\omega_{\widetilde\nu,A,\rho}$
is a basic form. 
On $P\times L_{\epsilon_0}$ the diffeomorphism
$\widetilde\psi_\nu$ is given by
$1_P\times U$, where $U$ lies in the
center of $U(k)$. Thus  
$ (1\times U)^*  \widetilde\omega_{A,\rho}=  \widetilde\omega_{A,\rho}$
and $\widetilde\psi_N$ is a symplectic diffeomorphism
of $(\widetilde M_N,\widetilde\omega_{N,\rho}).$

\end{proof}

Next we prove that the lift of a Hamiltonian diffeomorphism, under some constraints, is again a 
Hamiltonian diffeomorphism on the blown up manifold.
To that end, consider $\psi$  a Hamiltonian diffeomorphism on $(M,\omega)$
such that there exists a path $\{\psi_t\}_{0\leq t\leq 1}$ that start at the identity,
ends at $\psi$ and  for each $t$ the diffeomorphism $\psi_t$ is $N$-liftable.
Hence,  $\{\psi_{\nu,t}\}$  is a Hamiltonian path on 
$ (\mathcal{U}_{\epsilon_0}(\nu_0),\omega_{\nu,A})$. 
Moreover for each $t\in [0,1]$
there exist a matrix $U_t$ in the center of $U(k)$ 
such that 
$$
\psi_{\nu,t}[p,z]=[p,U_t(z)].
$$
Since each matrix $U_t$ is of the form $\lambda_t \cdot 1_{k\times k}$ for 
$\lambda_t\in \mathbb{C}$ and $|\lambda_t |=1$, it follows that the induced
Hamiltonian function of the path of unitary matrices is 
autonomous.
Moreover, it is given by 
$H(z):=c \sum_j  (x_j^2+y_j^2)=c| z|^2$ for some $c\in\mathbb{R}$.
Recall from Section \ref{ss:inducedHam}  that such Hamiltonian function
$H$ induces the Hamiltonian
$\widetilde H:L_{\epsilon_0}\to \mathbb{R}$ defined in Eq. (\ref{e:hamblow}).
Further it also satisfies the relation $\widetilde H(U(z),U(\ell))=\widetilde H(z,\ell)$.
Henceforth, we define the function
$\widetilde H_{t,N}:(\widetilde M_N,\widetilde\omega_{N,\rho}) \to \mathbb{R}$ by
\begin{eqnarray}
\label{e:hamblowN}
\widetilde H_{N,t}(q) :=
\left\{
	\begin{array}{ll}
		H_t(\widetilde\pi_{(M,N)}(q))   &  \textup{if }q\in 
					M\setminus\Phi(\overline{ \mathcal{U}_{\epsilon^\prime}(\nu_0)  })  \\
	\widetilde H(z,\ell) &   \textup{if }  q=[p,(z,\ell)] \in\mathcal{U}_{\epsilon_0}(\widetilde \pi^{-1}(\nu_0)). 
	\end{array}
\right.
\end{eqnarray}

Next we show that $\widetilde H_{N,t}$
 is the Hamiltonian function of the lifted path of symplectic diffeomorphisms.

\begin{prop}
Let  $\psi_t: (M,\omega)\to (M,\omega)$ be a path of Hamiltonian diffeomorphisms
such that each $\psi_t$ is $N$-liftable and $\psi_0=1$.
Then $\widetilde \psi_{N,t}: (\widetilde M_N,\widetilde\omega_{N,\rho})\to 
 (\widetilde M_N,\widetilde\omega_{N,\rho})$ is a Hamiltonian path
induced by $\widetilde H_{N,t}$. 
\end{prop}
\begin{proof}
Since the Hamiltonian path $\{ \psi_t \} $ is $N$-liftable, by Prop. \ref{p:liftsymplectic}
the path  $\{\widetilde\psi_{N,t}\}$ consists of symplectic diffeomorphisms on 
$(\widetilde M_N,\widetilde\omega_{N,\rho})$ and satisfies 
   $\psi_t\circ\widetilde\pi_{(M,N)}=\widetilde\pi_{(M,N)}\circ \widetilde\psi_{N,t}$ 
for each $t$. Thus if $\widetilde X_t$ is the time-dependent vector field
induced by $\{\widetilde\psi_{N,t}\}$, then  $(\widetilde\pi_{(M,N)})_*\widetilde X_t=X_t$.
Thus by the definition of the function $\widetilde H_{N,t}$ on (\ref{e:hamblowN}) and the fact that
$\widetilde\pi_{(M,N)}$ is a symplectic diffeomorphism on 
$M\setminus\Phi(\overline{ \mathcal{U}_{\epsilon^\prime}(\nu_0)  }) $
it follows that
$$
\widetilde\omega_N(\widetilde X_t ,\cdot)=d\widetilde H_{t,N}
$$
holds on $\widetilde M_N	\setminus \mathcal{U}_{\epsilon_0}(\widetilde \pi^{-1}(\nu_0))$.

Now on $ \mathcal{U}_{\epsilon_0}(\widetilde \pi^{-1}(\nu_0))\simeq P\times_{U(k)}L_{\epsilon_0}$
the path   $\{\widetilde \psi_{N,t}\}$ is given by $[p,(z,\ell)]\mapsto [p,(U_t(z),U_t(\ell))]$.
Thus the time-dependent vector field is actually independent of time and takes the form
$0+\widetilde X$ where $\widetilde X$ is a vector field on $L_{\epsilon_0}$.
Also, on 
$ \mathcal{U}_{\epsilon_0}(\widetilde \pi^{-1}(\nu_0))$
the symplectic form $\widetilde\omega_N$ is equal to 
$\widetilde\omega_{\widetilde \nu, A,\rho}=\widetilde\omega_{A,\rho}+\widetilde\pi_N^*(\omega)$.
Therefore,
\begin{eqnarray*}
\widetilde\omega_{\widetilde \nu, A,\rho}  (\widetilde X_t,\cdot)
 &=& \widetilde\omega_{\widetilde \nu, A,\rho}(0+\widetilde X,\cdot)\\
&=& \widetilde\omega_{A,\rho}(0+\widetilde X,\cdot)+\widetilde\pi_N^*(\omega)(0+\widetilde X,\cdot) \\
&=& 
\widetilde\omega_{\rho,\epsilon}(\widetilde X+0,\cdot)-  d\langle \widetilde\mu_\rho, A\rangle(0+\widetilde X,\cdot ) 
+\widetilde\pi_N^*(\omega)(0+\widetilde X,\cdot) \\
&=& 
\widetilde\omega_{\rho,\epsilon}(\widetilde X+0,\cdot)
\end{eqnarray*}
By
\cite[Prop. 3.7]{pedroza-hamiltonian}, on $L_{\epsilon_0}$ we have that
$\widetilde\omega_{\rho,\epsilon}(\widetilde X,\cdot)=d\widetilde H$ 
where the function $\widetilde H$ was
defined in Eq. (\ref{e:hamblow}).
Thus the proposition
follows.
\end{proof}

Notice that if $\{\psi_t\}$ is a loop of Hamiltonian diffeomorphisms on $(M,\omega)$
that is $N$-liftable, then the lift  $\{\widetilde \psi_{N,t}\}$ is also loop of Hamiltonian diffeomorphisms. 
We are interested in understanding the 
the behavior of the loop $\{\widetilde \psi_{N,t}\}$ in the fundamental group of $\textup{Ham}(\widetilde M_N,
\widetilde\omega_{N,\rho}).$

\section{Weinstein's morphism}
\label{s:wein}

The period group $\mathcal{P}(M,\omega)$ of $(M,\omega)$
is defined as the image of the pairing $\langle\omega,\cdot \rangle\colon H_2(M;\mathbb{Z})
\to \mathbb{R}$. Then Weinstein's morphism \cite{weinstein-coho},
$$
\mathcal{A}:\pi_1 (\textup{Ham}(M,\omega)) \to \mathbb{R}/ \mathcal{P}( M,\omega)
$$
is defined via the action functional on a loop $\psi=[\{\psi_t\}]$  as
$$
\mathcal{A}(\psi)=
\left[
-\int_{D} u^*( \omega) +\int_{0}^1 H_t(\psi_t(x_0)) dt \right].
$$
Here $D$ is the unit closed disk and  $u:D\to M$ is a smooth function such that
$u(\partial D)$ agrees with  the loop $\{\psi_t(x_0)\}$, $x_0\in M$  a base point and 
$H_t$ is the 1-periodic Hamiltonian function
induced by  the Hamiltonian  loop $\{\psi_t\}$ subject to the normalized condition
$$
\int_M H_t \, \omega^n=0
$$
for every $t\in [0,1]$.

The aim is to relate  Weinstein's morphism   on $\pi_1(\textup{Ham}(M,\omega))$ and
$\pi_1(\textup{Ham} (\widetilde M_N,\widetilde\omega_{N,\rho})      )$ on a particular class of loops.
To that end, recall that we must consider a Hamiltonian diffeomorphisms  on $(M,\omega)$
that are $N$-liftable. As seen above, such Hamiltonians  are
 induced by Hamiltonians functions on $(\mathbb{C}^n,\omega)$ that on a neighborhood of
the origin take the form
 $c \sum_j  (x_j^2+y_j^2)=c| z|^2$ for some $c\in\mathbb{R}$.

\begin{lem} 
\label{l:hamfuncrel}
Let $H:M\to\mathbb{R}$ be a smooth function that is induced by 
a Hamiltonian diffeomorphism that is $N$-liftable, $\rho$ as before and
$\widetilde H_N:\widetilde M_N\to\mathbb{R}$ the induced function defined in 
(\ref{e:hamblowN}). Then
$$
\int_{\widetilde M_N} \widetilde H_N \, \widetilde\omega_{N,\rho}^n
=\int_{M}  H\, \omega^n
-\int_{\Phi(\mathcal{U}_{\rho}(\nu_0))}  H\, \omega^n.
$$
\end{lem}
\begin{proof}
First note that $\widetilde M_N$ is the union of the disjoint sets
$M\setminus\Phi(\overline{ \mathcal{U}_{\epsilon^\prime}(\nu_0)  })$
and  $\mathcal{U}_{\epsilon^\prime}  (\widetilde\pi^{-1}(\nu_0)) $. Recall that
$\epsilon^\prime$ is such that
 $\rho+\epsilon<\epsilon^\prime < \epsilon_0.$
Then from the definition 
of $\widetilde H_N$ in (\ref{e:hamblowN}) and of $\widetilde H$ in (\ref{e:hamblow}) it follows that 
\begin{eqnarray*}
\int_{\widetilde M_N} \widetilde H_N \,\widetilde\omega_{N,\rho}^n
&=&
\int_{M\setminus\Phi(\overline{ \mathcal{U}_{\epsilon^\prime}(\nu_0)  }) }
 \hspace{-.6cm} \widetilde H_N \,\widetilde\omega_{N,\rho}^n
+\int_{\mathcal{U}_{\epsilon^\prime  }
(\widetilde\pi^{-1}(\nu_0))} \hspace{-.6cm} \widetilde H_N \,\widetilde\omega_{N,\rho}^n\\
&=&
\int_{M\setminus\Phi(\overline{ \mathcal{U}_{\epsilon^\prime}(\nu_0)  }) }
 \hspace{-.6cm} H\circ \widetilde\pi_{(M,N)} \,\widetilde\omega_{N,\rho}^n
+\int_{\mathcal{U}_{\epsilon^\prime  }
(\widetilde\pi^{-1}(\nu_0))} \hspace{-.6cm} \widetilde H \,\widetilde\omega_{N,\rho}^n\\
&=&
\int_{M\setminus\Phi(\overline{ \mathcal{U}_{\epsilon^\prime}(\nu_0)  }) }
 \hspace{-.6cm} H \,\omega^n
+\int_{\mathcal{U}_{\epsilon^\prime  }
(\widetilde\pi^{-1}(\nu_0))} \hspace{-.6cm} \widetilde H \,\widetilde\omega_{\widetilde\nu,A,\rho}^n\\
&=&
\int_{M\setminus\Phi(\overline{ \mathcal{U}_{\epsilon^\prime}(\nu_0)  }) }
 \hspace{-.6cm} H \,\omega^n
+\int_{(P\times_{U(k)}L_{\epsilon^\prime})\setminus \nu_0}
 \hspace{-.6cm}   \widetilde H \,\widetilde\omega_{\widetilde\nu,A,\rho}^n\\
\end{eqnarray*}

Recall that the diffeomorphism
 $F_{\rho,\epsilon}: \mathbb{C}^k\setminus \{0\}\to \mathbb{C}^k\setminus \overline{B}_\rho$
defined in 
(\ref{e:strech}) is $U(k)$-equivariant.
Hence it induces a diffeomorphism  
${\bf F}_{\rho,\epsilon}: P\times_{U(k)}(\mathbb{C}^k\setminus \{0\})
\to P\times_{U(k)}(\mathbb{C}^k\setminus \overline{B}_\rho)$ and
\begin{eqnarray*}
{\bf F}_{\rho,\epsilon} ^*(\omega_{\nu,A})  
&=&
{\bf F}_{\rho,\epsilon} ^*(\omega_A)  + {\bf F}_{\rho,\epsilon} ^*(\pi_N^*(\omega))\\
&=&
 \widetilde\omega_{A,\rho} + \widetilde\pi_N^*(\omega)\\
 &=&
\widetilde\omega_{\widetilde\nu,A,\rho}.
\end{eqnarray*}
It then follows that
\begin{eqnarray*}
\int_{(P\times_{U(k)}L_{\epsilon^\prime})\setminus \nu_0}
\hspace{-.6cm} \widetilde H \,\widetilde\omega_{\widetilde\nu,A,\rho}^n
&=&
\int_{\mathcal{U}_{\epsilon^\prime  }
(\widetilde\pi^{-1}(\nu_0))} \hspace{-.6cm} \widetilde H \, 
{\bf F}_{\rho,\epsilon} ^*\omega_{\nu,A}^n \\
&=&
\int_{P\times_{U(k)}(B_{\epsilon^\prime}\setminus \overline{B_\rho})}
 \hspace{-.6cm}  H \, 
\omega_{\nu,A}^n. 
\end{eqnarray*}


Hence
\begin{eqnarray*}
\int_{\widetilde M_N} \widetilde H_N \,\widetilde\omega_{N,\rho}^n
&=&
\int_{M\setminus\Phi(\overline{ \mathcal{U}_{\epsilon^\prime}(\nu_0)  }) }
 \hspace{-.6cm} H \,\omega^n
+\int_{P\times_{U(k)}(B_{\epsilon^\prime}\setminus \overline{B_\rho})}
 \hspace{-.6cm} H   \,\omega_{\nu,A}^n\\
 &=&
 \int_{M}  H\, \omega^n
-\int_{\Phi(\mathcal{U}_{\rho}(\nu_0))}  H\, \omega^n.
\end{eqnarray*}
\end{proof}

\section{Proof of Theorems
 \ref{t:mainrank} and \ref{t:weinsrelation}}
 \label{s:main}

Now we prove the main results of this article that appeared at the Introduction.

\smallskip
\noindent {\bf Theorem \ref{t:weinsrelation}.} {\em
Let $(M,\omega)$ be a closed symplectic manifold and $N\subset (M,\omega)$
a symplectic submanifold of codimension $2k>2$.
If $\{\psi_t\}$ is a loop of Hamiltonian diffeomorphisms on $(M,\omega)$
that is $N$-liftable to a loop $\{\widetilde\psi_{N,t}\}$
on $(\widetilde M_N,\widetilde \omega_{N,\rho})$, then
\begin{eqnarray*}
\mathcal{A}([\widetilde\psi_{N,t}]) = 
\left[\mathcal{A}([\psi_t]) +
\frac{1}{\textup{Vol}(\widetilde M_N,(\widetilde\omega_{N,\rho})^n)}
\int_0^1 
\int_{\Phi(\mathcal{U}_{\rho}(\nu_0))}  H_t\, \omega^n dt
 \right] 
\end{eqnarray*}
in $\mathbb{R}/\mathcal{P}(\widetilde M_N,\widetilde \omega_{N,\rho})$, where
$H_t$ if the normalized Hamiltonian function of the loop $\{\psi_t\}$. }

\begin{proof}
Let $\{\psi_t\}$ be a loop of Hamiltonian diffeomorphisms on $(M,\omega)$ that is
$N$-liftable with normalized Hamiltonian function $H_t$. If $\{\widetilde\psi_{N,t}\}$
is the induced loop on $(\widetilde M_N ,\widetilde\omega_{N,\rho})$, denote by
 $\widetilde H_{N,t}$ its Hamiltonian function defined by Eq. (\ref{e:hamblowN}).

It follows from Lemma \ref{l:hamfuncrel} that for each $t$,
\begin{eqnarray}
\label{e:intloc}
\int_{\widetilde M_N} \widetilde H_{N,t} \,\widetilde\omega_{N,\rho}^n
=
-\int_{\Phi(\mathcal{U}_{\rho}(\nu_0))}  H_t\, \omega^n.
\end{eqnarray}
Denote such number by $c(\rho,t)\in\mathbb{R}$. Hence the normalized Hamiltonian
function of the loop $\{\widetilde\psi_{N,t}\}$ is
$$
\widetilde{\bf H}_{N,t}:=\widetilde H_{N,t} -\frac{1}{\textup{Vol}(\widetilde M_N,
\widetilde\omega_{N,\rho}^n)} c(\rho,t).
$$

Thus, taking a disk $D$ and a point $x_0$ in $\widetilde M_N	\setminus \mathcal{U}_{\epsilon_0}(\widetilde \pi^{-1}(\nu_0))$,
\begin{eqnarray*}
\mathcal{A}([\widetilde\psi_{N,t}]) 
&=& 
\left[ -\int_D u^*(\widetilde\omega_{N,\rho}) + \int_0^1
\widetilde{\bf H}_{N,t} (\widetilde\psi_{N,t}(x_0)) dt\right]\\
 &=& 
\left[ -\int_D u^*(\omega) + \int_0^1
\widetilde H_{N,t} (\widetilde\psi_{N,t}(x_0)) dt \right. \\
& &\left. - \frac{1}{\textup{Vol}(\widetilde M_N,
\widetilde\omega_{N,\rho}^n)}    \int_0^1  c(\rho,t) dt \right].
\end{eqnarray*} 

Hence from Eq. (\ref{e:intloc}) we get
\begin{eqnarray*}
\mathcal{A}([\widetilde\psi_{N,t}]) = 
\left[\mathcal{A}([\psi_t]) +
\frac{1}{\textup{Vol}(\widetilde M_N,\widetilde\omega_{N,\rho}^n)}
\int_0^1 
\int_{\Phi(\mathcal{U}_{\rho}(\nu_0))}  H_t\, \omega^n dt
 \right].
\end{eqnarray*}
\end{proof}

In order to give meaning to the expression that was proven  above, now we relate the periods
groups of $(\widetilde M_N,\widetilde\omega_{N,\rho})$ and $( M,\omega)$.
Recall that the blow up map $\widetilde\pi_{(M,N)} : \widetilde M_N\to M $ is a diffeomorphism
on the complement of $\widetilde\pi_{(M,N)}^{-1}(N)$. Moreover,  $\widetilde\pi_{(M,N)}$
induces a symplectic diffeomorphism on the complement of a neighborhood of  
$\widetilde\pi_{(M,N)}^{-1}(N)$. Therefore, the maps
$$
\langle \widetilde\omega_{N,\rho}, \cdot \rangle: H_2(\widetilde M_N \setminus \widetilde\pi_{(M,N)}^{-1}(N);\mathbb{Z})
\to \mathbb{R} 
\textup{ and }
\langle \omega, \cdot \rangle: H_2( M_N;\mathbb{Z})
\to \mathbb{R} 
$$
have the same image. It remains to determine the value of 
$\langle \widetilde\omega_{N,\rho}, \cdot \rangle$ on the kernel of $\widetilde\pi_{(M,N)*}:
H_2(\widetilde M_N;\mathbb{Z})\to H_2( M;\mathbb{Z})$.
From the  commutative diagram
{\vskip .4cm}
$$
\xymatrix{
  & H_1(\nu\setminus\nu_0;\mathbb{Z}) \ar[d]_1\ar[l] &  H_2(\widetilde M_N;\mathbb{Z}) \ar[l] \ar[d]_{\widetilde\pi_{(M,N)*}}
&  H_2(M\setminus N;\mathbb{Z})\oplus H_2(\widetilde\nu;\mathbb{Z})\ar[d]_{1\oplus \widetilde\pi_*}\ar[l] & \ar[l]    \\
    & H_1(\nu\setminus\nu_0;\mathbb{Z}) \ar[l]        &  H_2(M;\mathbb{Z}) \ar[l]
&  H_2(M\setminus N;\mathbb{Z})\oplus H_2(\nu;\mathbb{Z}) \ar[l] & \ar[l]
}
$$
it suffices to look at the map $\widetilde\pi_* : H_2(\widetilde\nu;\mathbb{Z})
\to H_2(\nu;\mathbb{Z})$.

Recall that the closed 2-form $\widetilde\omega_{\widetilde\nu,A,\rho}$ on the
total space of the fibration    $\widetilde\pi_N:\widetilde\nu\to N$ 
restricts to  $\widetilde\omega_{\rho,\epsilon}$ on each fibre.

The relation between the homology of $\widetilde{M}_N$ and $M$ is as follows
$$
H_2(\widetilde M_N,\mathbb{Z})\simeq H_2(M,\mathbb{Z})\oplus \textup{ker}\{
{\widetilde\pi}_{N,*} :H_2(\widetilde \nu,\mathbb{Z})\to H_2(N,\mathbb{Z})\}.
$$
Moreover, if $\alpha\in H_2(\widetilde M_N;\mathbb{Z})$ takes the form
$\widetilde\pi_{(M,N)*}(\alpha)+ 0$ in the above decomposition, then
$$
\langle \widetilde\omega_{N,\rho},  \alpha \rangle =
\langle \omega,  \widetilde\pi_{(M,N)*}(\alpha)   \rangle.
$$
Since $\widetilde\pi_{(M,N)*}:
H_2(\widetilde M_N;\mathbb{Z})\to H_2( M;\mathbb{Z})$
is surjective, we have that
$$
\mathcal{P}( M,\omega)\subset \mathcal{P}(\widetilde M_N,\widetilde\omega_{N,\rho}).
$$
Finally, 
$\textup{ker}\{
{\widetilde\pi}_{N,*} :H_2(\widetilde \nu;\mathbb{Z})\to H_2(N;\mathbb{Z})\}$
is generated by a single class that its restriction to the fiber
is the generator of $H_2(\mathbb{C}P^{k-1};\mathbb{Z})$.
Hence we have the following result.

\begin{lem}
\label{l:relperiod}
Let $(M,\omega)$ and $(\widetilde M_N,\widetilde\omega_{N,\rho})$ as above.
Assume that $N$ is compact,
then
\begin{eqnarray*}
\mathcal{P}(\widetilde M_N,\widetilde\omega_{N,\rho})=
\mathcal{P}(M,\omega)+
\mathbb{Z}\langle\pi \rho^2 \rangle.
\end{eqnarray*}
\end{lem}

Recall  from Lemma \ref{l:sympformnormal} that on a neighborhood 
$P\times_{U(k)} B_{\epsilon_0}$
of the zero-section of
$\nu=P\times_{U(k)}\mathbb{C}^k$ the symplectic form is given by
\begin{eqnarray}
\label{e:expan}
\omega_{\nu,A}(u_1+v_1,u_2+v_2 ) &=&
\omega_0(v_1+X_{A(u_1)},v_2+X_{A(u_2)})- \label{omega0} \\
& & d\langle \mu, A\rangle(u_1+v_1,u_2+v_2 )+ \nonumber\\
& & \pi_N^*(\omega)(u_1+v_1,u_2+v_2 ) \nonumber
\end{eqnarray}
for $u_j\in TP$ and $v_j\in T\mathbb{C}^k$. Hence on the normal bundle $\nu$ we have that 
\begin{eqnarray}
\omega_{\nu,A}^n &=&   \sum_{j=0}^n 
\omega_0^j   \; \pi_N^*(\omega)^{n-j}   + d\alpha 
\end{eqnarray}
for some $(2n-1)$-form $\alpha$. (Keep in mind that in this equation  the 2-form
$\omega_0$ must be evaluated as in Eq. (\ref{e:expan})).   

\begin{lem}
\label{l:integralfunction}
Let $H:\nu\to \mathbb{R}$ be a smooth function. 
Then,
\begin{eqnarray*}
\int_{  \nu } H\; \omega_{\nu,A}^n 
&=& \int_{N}      (\pi_N)_*(   H\; \omega_0^k ) \;\omega^{n-k}.
\end{eqnarray*}
\end{lem}
\begin{proof}
First we integrate along the fiber with respect to $\pi_N: \nu\to N$ and then
use the commutativity of the diagram
$$
\xymatrix{
\nu \ar[r]^{\pi_N}
\ar[dr]&  N\ar[d]\\
     & \{pt\}.
}$$
Therefore by the projection formula we get that
\begin{eqnarray*}
(\pi_N)_*( H\; \omega_{\nu,A}^n ) &=&   \sum_{j=0}^n (\pi_N)_*\left(  
   H\; \omega_0^j   \; \pi_N^*(\omega)^{n-j}  \right )\\
&=&  \sum_{j=0}^n (\pi_N)_*(  
   H\; \omega_0^j ) \;\omega^{n-j}. \\
\end{eqnarray*}
Note that in this expression, the symplectic form $\omega$ is restricted to the submanifold $N$. 

Next we integrate on the submanifold $N$,
\begin{eqnarray*}
\int_{  \nu } H\; \omega_{\nu,A}^n 
&=& \sum_{j=0}^n  \int_{N}      (\pi_N)_*(  
   H\; \omega_0^j ) \;\omega^{n-j} \\
&=& \sum_{j=k}^n  \int_{N}      (\pi_N)_*(  
   H\; \omega_0^j ) \;\omega^{n-j} \\
&=& \int_{N}      (\pi_N)_*(  
   H\; \omega_0^k ) \;\omega^{n-k}.
\end{eqnarray*}
The last equality follows from the fact that the form $\omega_0$ comes from the 
fiber $\mathbb{C}^k$ of $\pi_N$.
\end{proof}

Notice that  from this result the evaluation of $\omega_0$ as in
Eq. (\ref{e:expan}) becomes irrelevant.

In particular, restrict the above integral to  ${\mathcal{U}_{\rho}(\nu_0)}$
and take
$H\equiv 1$. Then are able to compute the volume of ${\mathcal{U}_{\rho}(\nu_0)}$,
\begin{eqnarray*}
\textup{Vol}({\mathcal{U}_{\rho}(\nu_0)} , \omega_{\nu,A}^n )&=&
\int_{  {\mathcal{U}_{\rho}(\nu_0)} }  \omega_{\nu,A}^n \\
&=& \int_{N}      (\pi_N|_{\mathcal{U}_{\rho}(\nu_0)})_*(  \omega_0^k ) \;\omega^{n-k}\\
&=& \textup{Vol}(B_\rho , \omega_0^k )\cdot \textup{Vol}(N , \omega^{n-k} )\\
&=& \frac{\pi^{k} \rho^{2k}}{k!} \cdot \textup{Vol}(N , \omega^{n-k} ).
\end{eqnarray*}


\subsection{Description of the loop.}

In order to prove Thm. \ref{t:mainrank}, we need to defined a loop of Hamiltonian diffeomorphisms
on $(M,\omega)$ supported on $P\times_{U(k)}B_{\epsilon_0} $, a neighborhood of $N$. 
The definition of such loop follows the same line of thought as 
the Hamiltonian loop constructed in \cite[Sec. 2]{pea-rankham}.

Let $\alpha:\mathbb{R}\to\mathbb{R}$ such that $\alpha(0)=0$ and $\alpha(1)\neq 0$. Then consider
the path of diagonal  matrices $A_t:=\exp{(2\pi i\alpha(t))}\cdot 1_{k\times k}$. Hence $A_t$ lies
in the center of $U(k)$ and it induces on  $(\mathbb{C}^k,\omega_0)$
 a path $\{\psi^\alpha_t\}_{0\leq t\leq 1}$ of Hamiltonian diffeomorphisms
starting at the identity and
  with Hamiltonian function
$$
H^\alpha_t(z_1,\ldots , z_k):= \pi\; \alpha^\prime(t)\; \sum_{j=1}^k |z_j |^2.
$$
As before we have the positive  parameters $\rho,\epsilon_0$ and $\epsilon\in (0,1)$ such that
$\rho+\epsilon<\epsilon_0$, that are used in the definition of $(\widetilde M_N,\widetilde 
\omega_{N,\rho})$. 
Now fix a $U(k)$-invariant smooth function
 $g:\mathbb{C}^k\to \mathbb{R}$ supported in $B_{\epsilon_0}$ such that 
 $g\equiv 1$ on  $B_{\rho}$.
Consider the Hamiltonian function 
$$
H^{\alpha,g}_t: = g\cdot H^{\alpha}_t
$$
and let $\{\psi^{\alpha,g}_t\}_{0\leq t\leq 1}$ be the induced Hamiltonian path. Since $\omega_0$
and $H^{\alpha,g}_t$ are $U(k)$-equivariant, it follows that  each $\psi^{\alpha,g}_t$
is also $U(k)$-equivariant. Furthermore, on $B_\rho$ we have that 
$\psi^{\alpha,g}_t \equiv \psi^{\alpha}_t$ and $H^{\alpha,g}_t\equiv
H^{\alpha}_t$.
Note that as in \cite[Lemma 2.4]{pea-rankham},
\begin{eqnarray*}
\int_0^1 \int_{B_\rho}  \;H^{\alpha,g}_t\;\omega_0^k\;dt  
&= & (\alpha(1)-\alpha(0)) \;\pi
\int_{B_\rho} \sum_{j=1}^k x_j^2 + y_j^2 \frac{1}{k!}dx_1\cdots dy_k \\
&=&\alpha(1) \cdot \frac{\pi}{k!} \cdot \frac{\pi^{k} \rho^{2k+2}}{k+1}.
\end{eqnarray*}

Since each  $\psi^{\alpha,g}_t$ is $U(k)$-equivariant, it induces
a path of diffeomorphisms $1\times \psi^{\alpha,g}_t$ on $\nu=P\times_{U(k)}
\mathbb{C}^k.$ Furthermore, each diffeomorphism is
supported on $\mathcal{U}_{\epsilon}(\nu_0)=P\times_{U(k)} B_{\epsilon_0}$.

\begin{lem}
\label{l:hamnormal}
The path $1\times \psi^{\alpha,g}_t$ restricted to $\mathcal{U}_{\epsilon_0}(\nu_0)$
 is a  Hamiltonian path with Hamiltonian function $1\times H^{\alpha,g}_t.$
\end{lem}
\begin{proof}
Recall that the 2-form $\omega_{\nu,A}$ is 
a basic form on $P\times B_{\epsilon_0}$, in particular is $U(k)$-invariant.
Therefore, $1\times \psi^{\alpha,g}_t$ preserves the form $\omega_{\nu,A}$
 and is a symplectic diffeomorphism
on  $\mathcal{U}_{\epsilon_0}(\nu_0)$.

It remains to prove that is Hamiltonian. Note that by definition on $(B_{\epsilon_0},
\omega_0)$  we have that $\{\psi^{\alpha,g}_t\}$ is Hamiltonian with Hamiltonian 
function $H^{\alpha,g}_t$. Let $X^{\alpha,g}_t$ be the associated Hamiltonian vector
field, $\omega_0( X^{\alpha,g}_t, \cdot ) = dH^{\alpha,g}_t$. Recall from Eq.
(\ref{e:expan}) that 
\begin{eqnarray*}
\omega_{\nu,A}(u_1+v_1,u_2+v_2 ) &=&
\omega_0(v_1+X_{A(u_1)},v_2+X_{A(u_2)})- \label{omega0} \\
& & d\langle \mu, A\rangle(u_1+v_1,u_2+v_2 )+ \nonumber\\
& & \pi_N^*(\omega)(u_1+v_1,u_2+v_2 ) \nonumber
\end{eqnarray*}
for $u_j\in TP$ and $v_j\in TB_{\epsilon_0}$.  
If $X_{\xi}$ is the vector field on $B_{\epsilon_0}$ induced
by $\xi\in \mathfrak{u}(k)$ we have that
$\omega_0(X_{\xi_1},X_{\xi_2})=\langle \mu , [{\xi_1}, {\xi_2}]\rangle $
for any $\xi_1,\xi_2\in \mathfrak{u}(k)$. In particular, the vector field $X^{\alpha,g}_t$
is of the form $X_{\xi}$ for ${\xi}$ in the center of $ \mathfrak{u}(k)$ and therefore
$\omega_0(X^{\alpha,g}_t ,X_{\xi})=0$ for any ${\xi}$. Then for 
$u\in TP$ and   $v\in TB_{\epsilon_0} $ 
\begin{eqnarray*}
\omega_{\nu,A}( X^{\alpha,g}_t,  u+v ) 
&=&\omega_0(X^{\alpha,g}_t ,v+X_{A(u)})  \\
&=&\omega_0(X^{\alpha,g}_t ,v)  \\
&=&  dH^{\alpha,g}_t  (    v  ) \\
&=& d(1\times H^{\alpha,g}_t)  (u,v)
\end{eqnarray*}
where $\langle \mu, A\rangle$ vanishes since $v_1=0$.
\end{proof}

Now we will consider another path of Hamiltonian diffeomorphisms in order
to create a loop by concatenating it with the path 
$\{\psi^{\alpha,g}_t\}_{0\leq  t\leq 1}$. Hence, consider  $\beta:\mathbb{R}\to\mathbb{R}$ 
such that $\beta(0)=0$ and  $\alpha(t)-\beta(t)\equiv1$ in a neighborhood of 1.
As above denote
by  $\{\psi^{\beta,g}_t\}_{0\leq t\leq 1}$ the induced the path of Hamiltonian diffeomorphisms
 with Hamiltonian function $
H^{\beta,g}_t$ supported on $B_{\epsilon_0}$. Then the desired loop of Hamiltonian
 diffeomorphisms  $\{  
\psi^{\mathbb{C}^k}_t  \}_{0\leq t\leq 2}$    is defined as
\begin{eqnarray}
\label{e:loop}
\psi^{\mathbb{C}^k}_t :=
\left\{
	\begin{array}{ll}
		\psi^{\alpha,g}_t  & t\in[0,1] \\
		\psi^{\beta,g}_{2-t}  & t\in[1,2]
	\end{array}
\right.
\end{eqnarray}
with Hamiltonian function 
\begin{eqnarray}
\label{e:Hloop}
H^{\mathbb{C}^k}_t :=
\left\{
	\begin{array}{ll}
		H^{\alpha,g}_t  & t\in[0,1] \\
		H^{\beta,g}_{2-t}  & t\in[1,2].
	\end{array}
\right.
\end{eqnarray}

Furthermore, we have that 
\begin{eqnarray*}
\int_0^2 \int_{B_\rho}  \; H^{\mathbb{C}^k}_t   \;\omega_0^k\;dt  
= (\alpha(1)-\beta(1)) \; \frac{\pi^{k+1} \rho^{2k+2}}{(k+1)!}
= \frac{  (\pi \rho^2)^{k+1}}{(k+1)!}.
\end{eqnarray*}

We have all the prerequisites to proceed with the proof of the main theorem.

\medskip
\noindent {\bf Theorem \ref{t:mainrank}.} {\em
Let $(M,\omega)$ be a rational closed symplectic manifold and $N\subset (M,\omega)$
a closed symplectic submanifold of codimension $2k>2$.
Let $(\widetilde M,\widetilde\omega_{N,\rho})$ be the symplectic blow up
of $(M,\omega)$ along $N$. Then for some small values of $\rho$ the rank of 
$\pi_1(\textup{Ham}(\widetilde M,\widetilde\omega_{N,\rho}))$
is positive.}

\begin{proof}
Consider the loop  of Hamiltonian diffeomorphisms $\{  \psi^{\mathbb{C}^k}_t  \}_{0\leq t\leq 2}$
defined above.
Since each  $\psi^{\mathbb{C}^k}_t$ is $U(k)$-equivariant, it induces
a loop of diffeomorphisms $\{\psi^{\nu}_t\}$ on $\nu=P\times_{U(k)}
\mathbb{C}^k.$ Furthermore, each diffeomorphism is
supported on $P\times_{U(k)} B_{\epsilon_0}$.

From Lemma \ref{l:hamnormal}, we have that $\{1\times \psi^{\nu}_t\}$
is a loop of Hamiltonian diffeomorphisms on $P\times_{U(k)} B_{\epsilon_0}$
with Hamiltonian function $H^\nu_t=1\times H^{\mathbb{C}^k}_t.$ 
Since $P\times_{U(k)} B_{\epsilon_0}$ is 
symplectomorphic to a neighborhood $\mathcal{U}(N)$, we consider 
$\{\psi^{\nu}_t\}$ as a loop of Hamiltonian diffeomorphisms on $(M,\omega)$.
Denote such loop by $\{\psi^M_t\}$ with Hamiltonian function $H^M_t$.
Notice that such loop is contractible in $\textup{Ham}(M,\omega)$ and the
Hamiltonian function is not normalized. Thus consider the function
\begin{eqnarray}
\label{e:normalconstant}
c_t:= \int_M  H^M_t\; \omega^n
= \int_{\mathcal{U}(N)} H^M_t\; \omega^n.
\end{eqnarray}
Then the normalized Hamiltonian function of $\{\psi^M_t\}$ is defined as
\begin{eqnarray*}
H^{M,\textup{norm}}_t:=  H^M_t -  \frac{c_t}{\textup{Vol}(M,\omega^n)}.
\end{eqnarray*}

From its definition, the loop  $\{\psi^M_t\}_{0\leq t\leq 2}$ is $N$-liftable. 
Denote by $\{\widetilde\psi_{N,t}\}_{0\leq t\leq 2}$ the induced loop of
Hamiltonian diffeomorphisms.
Hence  by Thm \ref{t:weinsrelation} we have that
\begin{eqnarray*}
\mathcal{A}([\widetilde\psi_{N,t}]) = 
\left[
\frac{1}{\textup{Vol}(\widetilde M_N,\widetilde\omega_{N,\rho}^n)}
\int_0^2 
\int_{\Phi(\mathcal{U}_{\rho}(\nu_0))}  H^{M,\textup{norm}}_t\, \omega^n dt
 \right] 
\end{eqnarray*}
in $\mathcal{P}(M,\omega)+
\mathbb{Z}\langle\pi \rho^2 \rangle.$ 
From the definition of the symplectic blow up we know that 
$  \textup{Vol}(\widetilde M_N,\widetilde\omega_{N,\rho}^n)$ is equal to 
$\textup{Vol}(M,\omega^n) -
\textup{Vol}  ( \mathcal{U}_{\rho}(\nu_0), \omega_{\nu,A}^n )$.

Next we compute the integrals. To that end, by the definition
of the Hamiltonian function on the neighborhood 
$\Phi(\mathcal{U}_{\rho}(\nu_0))$
of $N$ we 
have that
\begin{eqnarray*}
\int_0^2 
\int_{\Phi(\mathcal{U}_{\rho}(\nu_0))}  H^{M,\textup{norm}}_t\, \omega^n \, dt
&=&
\int_0^2 
\int_{\Phi(\mathcal{U}_{\rho}(\nu_0))}   H^M_t -   \frac{c_t}{\textup{Vol}(M,\omega^n)}   \,   \omega^n  \,dt \\
&=&
\int_0^2 
\int_{\mathcal{U}_{\rho}(\nu_0)}       1\times H^{\mathbb{C}^k}_t     \,  \omega_{\nu,A}^n  \,dt    \\
 &  &   -          
    \frac{ \textup{Vol}  ( \mathcal{U}_{\rho}(\nu_0), \omega_{\nu,A}^n )  }{\textup{Vol}(M,\omega^n)} \int_0^2
c_t   \, dt.
\end{eqnarray*}
By replacing the normal bundle $\nu$ by the disk bundle $\mathcal{U}_{\rho}(\nu_0)$ on the 
first sumand, we get
by Lemma \ref{l:integralfunction} that
\begin{eqnarray*}
\int_{\mathcal{U}_{\rho}(\nu_0)}       1\times H^{\mathbb{C}^k}_t      \omega_{\nu,A}^n    
&=& 
\int_{N}      (\pi_N  |_{\mathcal{U}_{\rho}(\nu_0)} )_*(    1\times H^{\mathbb{C}^k}_t   \; \omega_0^k ) \;\omega^{n-k}.
\end{eqnarray*}
Since  
\begin{eqnarray*}
 (\pi_N  |_{\mathcal{U}_{\rho}(\nu_0)} )_*(    1\times H^{\mathbb{C}^k}_t   \; \omega_0^k )
&=&
 \int_{B_\rho}  \; H^{\mathbb{C}^k}_t   \;\omega_0^k,
\end{eqnarray*}
it follows that 
\begin{eqnarray*}
\int_0^2 \int_{\mathcal{U}_{\rho}(\nu_0)}       1\times H^{\mathbb{C}^k}_t      \omega_{\nu,A}^n  \,dt   
&=& 
 \frac{  (\pi \rho^2)^{k+1}}{(k+1)!}\textup{Vol}(N,\omega^{n-k}) 
\end{eqnarray*}
since $\alpha(1)-\beta(1)=1.$

Notice that $c_t$ can be rewritten as
\begin{eqnarray*}
 {c_t}
&=& 
\int_{M} H^M_t\; \omega^n\\\
&=& 
\int_{\mathcal{U}(N)} H^M_t\; \omega^n\\
&=& 
\int_{\nu}  1\times H^{\mathbb{C}^k}_t 
 \; \omega_{\nu,A}^n \\
&=& 
\int_{N}      (\pi_N   )_*(    1\times H^{\mathbb{C}^k}_t   \; \omega_0^k ) \;\omega^{n-k}.
 \end{eqnarray*}
Before the fiber was $(B_\rho,\omega_0)$, in this case the fiber is $(\mathbb{C}^k,\omega_0)$ 
and  $ H^{\mathbb{C}^k}_t $ has  compact support. 
But on  $(\mathbb{C}^k,\omega_0)$, we know that the 
Calabi morphism vanishes. Therefore we have that
$$
\int_0^2    {c_t} \; dt=
\int_0^2 (\pi_N   )_*(    1\times H^{\mathbb{C}^k}_t   \; \omega_0^k )\; dt =0
$$
and 
\begin{eqnarray*}
\mathcal{A}([\widetilde\psi_{N,t}]) = 
\left[
\frac{1}{\textup{Vol}(M,\omega^n) -
  \frac{\pi^{k} \rho^{2k}}{k!} \cdot \textup{Vol}(N , \omega^{n-k} )
 } \cdot
 \frac{  (\pi \rho^2)^{k+1}}{(k+1)!}\textup{Vol}(N,\omega^{n-k}) 
 \right]
\end{eqnarray*}
in $\mathcal{P}(M,\omega)+\mathbb{Z}\langle\pi \rho^2 \rangle.$ By hypothesis, 
$\textup{Vol}(M,\omega^n),  \textup{Vol}(N,\omega^{n-k})$ and
 $\mathcal{P}(M,\omega)$  are in $\mathbb{Q}$.
Therefore, it follows from the above  expression that if $\rho$ is such that $\pi\rho^2$ is a 
transcendental number
then $\mathcal{A}([\widetilde\psi_{N,t}]^m)\neq 0$ in 
$\mathcal{P}(M,\omega)+\mathbb{Z}\langle\pi \rho^2 \rangle$
for all non zero $m\in\mathbb{Z}$. Henceforth, $\pi_1(\textup{Ham}(\widetilde M,\widetilde\omega_{N,\rho}))$
has positive rank. 
\end{proof}


\bibliographystyle{acm}
\bibliography{/Users/andres/Dropbox/Documentostex/Ref.bib} 
\end{document}